\documentclass[letterpaper, 10 pt, conference]{ieeeconf}  

\IEEEoverridecommandlockouts                              

\usepackage{graphics} 
\usepackage{amsmath} 
\usepackage{amssymb}  
\usepackage{mathtools}
\usepackage{flushend}

\usepackage[dvipsnames]{xcolor}
\usepackage{pgfplots}
\usepackage{algpseudocode}
\usepackage{algorithm}
\pgfplotsset{compat=1.16}
\usepgfplotslibrary{groupplots}

\newcommand{\Id}{{\mathrm{Id}}}

\newcommand{\tom}[1]{\textcolor{black}{#1}}

\pgfplotscreateplotcyclelist{colors}{
        {CornflowerBlue},
        {BurntOrange},
        {RubineRed},
        {Emerald},
        {SpringGreen},
        {Goldenrod}
}

\usepackage{./TLC}

\title{\LARGE \bf
Circuit Analysis using Monotone+Skew Splitting
}

\author{Thomas Chaffey$^{1}$, Sebastian Banert$^{2}$, Pontus  Giselsson$^{2}$ and Richard Pates$^{2}$
\thanks{*S. Banert and P. Giselsson were financially supported by the Swedish Research Council. S. Banert, P. Giselsson, and R. Pates are members of ELLIIT.  R. Pates and T. Chaffey were financially supported by funding from ERC grant agreement No 834142.}
\thanks{$^{1}$T. Chaffey is with Pembroke College,
        University of Cambridge, United Kingdom, {\tt\small tlc37@cam.ac.uk}}%
\thanks{$^{2}$S. Banert, P. Giselsson and R. Pates are with the Department of Automatic Control, Lund University, Lund, Sweden
        {\tt\small \{sebastian.banert, pontusg, richard.pates\}@control.lth.se}}%
}

\begin{document}

\maketitle
\thispagestyle{empty}
\pagestyle{empty}

\begin{abstract}
It is shown that the behavior of an $m$-port circuit of maximal monotone elements can be expressed as a zero of the sum of a maximal monotone operator containing the circuit elements, and a structured skew-symmetric linear operator representing the interconnection structure, together with a linear output transformation.  The Condat--V\~u algorithm solves inclusion problems of this form, and may be used to solve for the periodic steady-state behavior, given a periodic excitation at each port, using an iteration in the space of periodic trajectories.
\end{abstract}

\section{INTRODUCTION}\label{sec:intro}

The study of monotone operators originated in attempts to generalize classical electrical network theory from circuits of linear, time invariant elements to circuits containing nonlinear elements.  Building on the ``quasi-linear resistor'' of Duffin \cite{Duffin1946}, Minty defined a monotone operator as an incremental, nonlinear generalization of a passive LTI circuit component \cite{Minty1960, Minty1961}.  The theory of monotone operators has since become central to several lines of research.  Following the influential work of Rockafellar \cite{Rockafellar1967, Rockafellar1976}, monotone operators have become a pillar of optimization theory, the essential observation being that the subdifferential of a closed, convex and proper function is a maximal monotone operator. \emph{Splitting algorithms} for solving convex optimization problems have seen a surge in interest in the last decade, due to their amenability to large-scale and distributed problems \cite{Rockafellar1997, Ryu2016, Ryu2022, Parikh2013, Bertsekas2011, Combettes2011}.  There is a large body of literature exploiting the property of monotonicity in the time-stepping simulation of discontinuous dynamical systems \cite{Brogliato2020, Camlibel2016, Brogliato2004, Brezis1973}, \tom{and a circuit interpretation of monotonicity has been used to study consensus \cite{Burger2014}}.

The use of monotone operators in circuit analysis has recently been revisited by Chaffey and Sepulchre \cite{Chaffey2021d}.  A one-port circuit consisting of a collection of monotone operators connected in series and parallel defines a monotone current--voltage relation.  A splitting algorithm is constructed which solves the steady-state behavior of the circuit, where the structure of the algorithm is in direct correspondence with the circuit topology.  This gives a method for solving the forced periodic response of a nonlinear circuit using an iteration in the space of periodic trajectories, rather than forwards in time.  

Alternative methods for finding the periodic response of nonlinear circuits are either approximate (for instance, the describing function method \cite{Krylov1947}), or involve integrating a differential equation forwards in time and waiting for convergence to steady state \cite{Aprille1972, Cellier2006}.  Specialized methods for circuits containing monotone elements may be found in the literature on non-smooth differential inclusions \cite{Acary2008}.  The methods of Heemels \emph{et al.} \cite{Heemels2017} involve iteratively computing the resolvent of the transition map of a differential inclusion, \tom{and most closely resemble the method presented in this paper}.  For circuits modelled as linear complementarity systems, methods have been developed by Ianelli \emph{et al.} \cite{Iannelli2011} and Meingast \emph{et al.} \cite{Meingast2014}.  

In this paper, we generalize the signal space approach of \cite{Chaffey2021d} to circuits with an arbitrary number of ports and an arbitrary interconnection structure.  For such a circuit, the mapping between any two voltages or currents is not necessarily monotone.  However, the circuit can be described by a monotone mapping from a set of external excitations to a set of internal currents and voltages, together with a linear mapping from the internal currents and voltages to a set of external responses.  We use a splitting algorithm to solve for the internal signals in terms of the excitations, and then apply the output mapping to obtain the responses.  The interconnection structure of the circuit is described by a skew-symmetric linear operator (\tom{the graph of which defines a Dirac structure} \cite[$\S 2.2.2$]{vanderSchaft2014}), and the collection of elements is described by a diagonal monotone operator.  Partitioning the elements into admittances (mapping voltages to currents) and impedances (mapping currents to voltages) expresses the circuit as precisely the \emph{monotone+skew form} \cite{Briceno-Arias2010} of zero inclusion which may be solved using the Condat--V\~u algorithm \cite{Condat2013, Vu2013}, a generalization of the Chambolle--Pock algorithm \cite{Chambolle2011}.
In fact the pairing between optimization problem and electrical circuit is stronger, and the solutions to any zero inclusion problem of this form define the behavior of an electrical circuit that can be synthesized with passive and reciprocal components.  \tom{To the best of the authors' knowledge, the method of \cite{Chaffey2021d} and the generalization we present in this paper represent the first use of splitting algorithms to explicitly partition a circuit simulation problem according to the interconnection structure.  This gives a method of circuit simulation which is readily distributed and comes pre-equipped with an algorithmic complexity theory.  A further advantage is that set-valued elements, such as the ideal diode, are handled via their resolvents, which are continuous functions.}

The remainder of this paper is structured as follows.  Section~\ref{sec:problem_statement} introduces the class of $m$-port circuits and the problem considered in this paper, solving the periodic behavior of an $m$-port circuit.  Section~\ref{sec:representation} describes how the behavior of an $m$-port circuit of monotone elements may be expressed as a \tom{monotone+skew} zero inclusion: this is captured in Theorem~\ref{thm:representation}, which represents the main theoretical contribution of this paper.  Section~\ref{sec:algorithms} gives an overview of the Condat--V\~u algorithm, which we use in Section~\ref{sec:example} to two detailed examples of the proposed method.

\section{PRELIMINARIES AND PROBLEM STATEMENT}\label{sec:problem_statement}

Let $\mathcal{H}$ denote an arbitrary Hilbert space, and $\ip{\cdot}{\cdot}$ denote its inner product.  We will focus in particular on the Hilbert space $L^m_{2, T}$ of square-integrable signals defined on the time axis $[0, T]$ and taking values in $\R^m$, with inner product
\begin{IEEEeqnarray*}{rCl}
    \ip{u}{y} &=& \int_0^T u(t)\tran y(t) \dd{t}.
\end{IEEEeqnarray*}
$T$-periodic signals may be represented as signals in $L^m_{2, T}$ by considering only a single period.  The reason for using the space $L_{2, T}$ is that it allows our method to be made computational via discretization, as demonstrated in Section~\ref{sec:example}.  The representation of a circuit introduced in Section~\ref{sec:representation} is equally valid on any Hilbert space.

By an \emph{operator} on a Hilbert space $\mathcal{H}$, we will mean a (possibly) multi-valued mapping on $\mathcal{H}$, that is, a function $\mathcal{H} \to 2^\mathcal{H}$, where $2^\mathcal{H}$ denotes the power set of $\mathcal{H}$.  Given an operator $A: \mathcal{H} \to 2^\mathcal{H}$, the \emph{relation} of $A$, denoted $\operatorname{rel}(A)$, is defined to be the set $\{(u, y) \; | \; y \in A(y)\} \subseteq \mathcal{H}\times\mathcal{H}$.  An operator is said to be \emph{monotone}\footnote{Monotonicity on $L_{2, T}$ for all $T$ is equivalent to incremental passivity \cite{Desoer1975}.} if, for any $(u_1, y_1), (u_2, y_2) \in \operatorname{rel}(A)$, 
\begin{equation}\label{eq:monotone}
\ip{u_1 - u_2}{y_1 - y_2} \geq 0.
\end{equation}  A monotone operator is said to be \emph{skew} or \emph{lossless} if \eqref{eq:monotone} holds with equality.  A monotone operator is said to be \emph{maximal} if its relation is not properly contained in the graph of any other monotone operator.  

The identity operator, $x \mapsto x$, is denoted $\Id$.  The \emph{(relational) inverse} of an operator $A$ is the operator $A^{-1}$ with relation $\operatorname{rel}(A^{-1}) = \{ (y, u) \; | \; (u, y) \in \operatorname{rel}(A)\}$.  The $\alpha$-\emph{resolvent} of an operator $A$ is the operator $(\Id + \alpha A)^{-1}$.

By an operator on $L^m_{2, T}$, we will mean an operator which maps a $T$-periodic signal into another $T$-periodic signal, both of which have a finite integral over $[0, T]$.  This mapping can be completely described by considering a single period of the input and a single period of the output, giving a mapping on $L^m_{2, T}$.

A \emph{one-port circuit} is a circuit with two terminals, described by a relation between the voltage between the two terminals, and the current through them.
We will model circuit elements as one-port circuits described by a maximal monotone operator on $L_{2, T}$, either from voltage to current (admittance form) or current to voltage (impedance form).  Particular examples include the admittance of a nonnegative, LTI capacitor, the impedance of a nonnegative, LTI inductor, and any maximal monotone nonlinear resistor \cite[$\S$VI]{Chaffey2021d}.  The relations for these elements are defined respectively through laws of the form
\begin{IEEEeqnarray*}{rCl}
    i(t) &=& \td{}{t} C v(t)\\
    v(t) &=& \td{}{t} L i(t)\\
    v(t) &=& R(i(t)),
\end{IEEEeqnarray*}
where $C \geq 0$, $L \geq 0$ and $R(\cdot)$ is non-decreasing and defines a continuous curve in the plane, with no endpoints.

An ideal $(m + n)$-port transformer is a device with a hybrid representation of the following form
\begin{IEEEeqnarray}{rCl}
    \begin{pmatrix} \mathbf{i}_1(t) \\ \mathbf{v}_2(t) \end{pmatrix} &=& \begin{pmatrix} 0 & T\tran \\ -T & 0 \end{pmatrix}\begin{pmatrix} \mathbf{v}_1(t) \\ \mathbf{i}_2(t) \end{pmatrix},\label{eq:transformer}
\end{IEEEeqnarray}
where $T \in \R^{n + m}$ is the \emph{turns ratio matrix}, and each entry \tom{$a_{pq}$ of $T$ gives the ratio of windings between driving point $p$ and driving point $q$.}

\begin{figure}
\centering
    \includegraphics{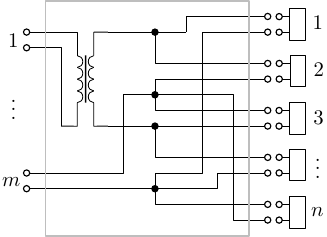}
    \caption{An $m$ port circuit modelled using the element extraction approach. The $m$ terminal pairs on the left represent external driving points.  The grey box contains wires and ideal transformers, governed by Kirchoff's laws and the ideal transformer relation \eqref{eq:transformer}.  The $n$ one-port elements are connected to $n$ internal ports, shown on the right.}
    \label{fig:element_extraction}
\end{figure}

We define a \emph{monotone $m$-port circuit with $n$ elements} to be a circuit constructed in the following way.  Begin with a box of wires and ideal transformers, connected arbitrarily.  From this box, draw $m$ pairs of terminals to form $m$ external \emph{driving points}.  Draw a further $n$ pairs of terminals, and across each of these pairs, connect a monotone one-port circuit element.  This structure is shown in Figure~\ref{fig:element_extraction}.
Throughout the paper, we will orient the voltage and current at a terminal pair of the box of wires so that $i(t)v(t)$ represents the instantaneous power \emph{extracted} from the box (this is opposite to the usual orientation, but means we can choose the usual orientation for the circuit elements when we perform element extraction in Section~\ref{sec:representation}).  An excitation at a driving point refers to an applied voltage (respectively, current), the response refers to the resulting observed current (voltage).

If we collect the driving point currents and current through each device in a vector $\iota \in L_{2, T}^{m + n}$ and the driving point voltages and voltage across each component into a vector $\nu \in L_{2, T}^{m + n}$, the \emph{behavior} \cite{Willems2007} of an $m$-port circuit $\mathcal{C}$ on $L_{2, T}$ is the set of all $\{(\iota, \nu) \in L_{2, T}^{m + n} \times L_{2, T}^{m + n}\}$ that obey the circuit interconnection constraints and device laws of $\mathcal{C}$.
Theorem~\ref{thm:representation} of Section~\ref{sec:representation} will give a precise mathematical formulation for this informal definition of a circuit and its behavior.
The problem considered in the remainder of this paper is the following:

\begin{quote}\emph{
given a monotone $m$-port circuit $\mathcal{C}$, and a $T$-periodic excitation signal at each driving point, find a set of steady-state periodic responses, if such a set exists.
}\end{quote}

It will be shown in Theorem~\ref{thm:representation} that this problem is equivalent to finding a zero to a sum of monotone operators and applying a linear transformation (Equations~\eqref{eq:inclusion1} and~\eqref{eq:inclusion2}).

\tom{
We conclude this preliminary section with a simple example, which illustrates the circuit modelling approach formalized in Section~\ref{sec:representation}.
}

\begin{example}\label{ex:basic}
\tom{
   \begin{figure}
        \centering
        \includegraphics{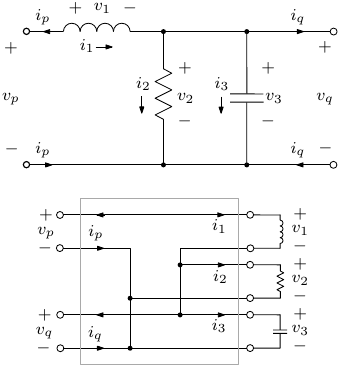}
        \caption{\tom{Above: a two-port RLC circuit. Below: the same circuit, with the elements extracted.}}
        \label{fig:basic_example}
    \end{figure}
    Consider the RLC circuit shown in the top of Figure~\ref{fig:basic_example}.  Assume all signals are $T$-periodic for some fixed $T > 0$, and belong to $L_{2, T}$.  Suppose we want to determine the mapping from $v_p$ to $v_q$, assuming that $i_q(t) = 0$ for all $t$ and $i_p$ is free.  The bottom of Figure~\ref{fig:basic_example} shows the same circuit, with the elements on one side, the external ports on the other, and a collection of wires joining them, enclosed in the grey dashed box.  Kirchoff's current and voltage laws for the grey box can be written as
    \begin{IEEEeqnarray}{rCl}
        \begin{pmatrix} i_p(t) \\ v_q(t) \\ v_1(t) \\ v_2(t) \\ i_3(t) \end{pmatrix}
        &{=}& \left(\begin{array}{c c | c c c} 
                0 & 0 & 0 & -1 & -1 \\
                0 & 0 & -1 & 0 & 0  \\\hline
                1 & 0 & 0 & 0 & -1  \\
                0 & 0 & 0 & 0 & 1  \\
                0 & -1 & 1 & -1 & 0  
\end{array}\right)\begin{pmatrix} v_p(t) \\ i_q(t) \\ i_1(t) \\ i_2(t) \\ v_3(t) \end{pmatrix}.\label{eq:ex1_KCL}
    \end{IEEEeqnarray}
    This is a \emph{hybrid representation} for the subcircuit contained in the grey box.  A partition of the circuit variables which admit a hybrid representation is guaranteed to exist, but in general is not unique.  The standard methods of loop and cut-set analysis can be used to find a hybrid representation \cite{Desoer1969}.  Here, we have put the variables we would like to solve, $v_q$ and $i_p$, on the left, in the \emph{response} vector, and the variables we would like to treat as inputs, $v_p$ and $i_q$, on the right, in the \emph{excitation} vector.
    }


\tom{
    Suppose $L, R, C \geq 0$ are the inductance, resistance and capacitance respectively, of the elements in Figure~\ref{fig:basic_example}.
    Let $\hat{L}$ and $\hat{R}$ denote the impedances of the inductor and resistor respectively, as operators on $L_{2, T}$.  Similarly, let $\hat{C}$ denote the admittance of the capacitor, as an operator on $L_{2,T}$.
    Substituting the device laws into Equation~\eqref{eq:ex1_KCL} then gives a complete description of the circuit behavior, as the sum of a monotone operator (offset by the inputs $v_p$ and $i_q$) and a skew operator, together with a linear output transformation:
    \begin{IEEEeqnarray*}{*rCr+x*}
        0 &\in& \begin{pmatrix} \hat{L}(i_1) \\ \hat{R}(i_2) \\ \hat{C}(v_3) \end{pmatrix}
        + \left(\begin{array}{c c c} 
                0 & 0 & -1  \\
                0 & 0 &  1  \\
                1 & -1 & 0  
    \end{array}\right)\otimes \Id \begin{pmatrix} i_1 \\ i_2 \\ v_3 \end{pmatrix}\IEEEyesnumber \label{eq:ex1_inclusion}\\
        &&+ \begin{pmatrix} 1 & 0 \\ 0 & 0 \\ 0 & -1 \end{pmatrix}\otimes \Id \begin{pmatrix}v_p \\ i_q \end{pmatrix},\\
        &&\begin{pmatrix} i_p \\ v_q \end{pmatrix} =
        \begin{pmatrix} 0 & -1 & -1 \\ -1 & 0 & 0 \end{pmatrix}\otimes \Id \begin{pmatrix} v_1 \\ i_2 \\ i_3 \end{pmatrix}, 
    \end{IEEEeqnarray*}
    where $\Id$ is the identity mapping on $L_{2,T}$.  
    }
\end{example}

\section{FROM CIRCUITS TO MONOTONE INCLUSIONS}\label{sec:representation}

In this section, we formalize the modelling approach of Example~\ref{ex:basic}.
We adopt the element extraction approach of Hughes \cite{Hughes2017c}, which builds on
the approach of Anderson and Newcomb \cite{Anderson1966}.  This allows us to express the behavior of an $m$-port circuit as the zeros of the sum of a skew-symmetric linear operator, representing the interconnection structure, and a diagonal monotone operator, representing the circuit elements, together with a linear mapping from the internal signals to the response signals.  
%
%
%
This representation is formalized in the following theorem.

\begin{theorem}\label{thm:representation}
    Let $\mathcal{C}$ be a monotone $m$-port circuit with $n$ elements.  Then there exist:
    \begin{enumerate}
        \item a partition of the driving point voltages $\{v_k\}_{k = 1\ldots m}$ and driving point currents $\{i_k\}_{k = 1 \ldots m}$ into vectors $\mathbf{u}, \mathbf{y} \in L_{2, T}^m$,
        \item vectors $\mathbf{i} \in L_{2, T}^{p}$ and $\mathbf{v} \in L_{2, T}^{q}$, $p + q = n$, such that for each element either the voltage across it is in $\mathbf{v}$ or the current through it is in $\mathbf{i}$,
        \item matrices $\tilde{M} \in \R^{q \times p}$, $\tilde{B}_R \in \R^{p \times m}$, $\tilde{B}_G \in \R^{q \times m}$ and $\tilde{D} \in \R^{m \times m}$, and
        \item monotone operators $R: L^p_{2, T} \to L^p_{2, T}$, and $G: L^q_{2, T} \to L^q_{2, T}$, which are the concatenation of the impedances of elements whose currents are in $\mathbf{i}$ and admittances of elements whose voltages are in $\mathbf{v}$,
    \end{enumerate}
    such that the behavior of $\mathcal{C}$ on $L_{2, T}$ is the set of solutions to
\begin{IEEEeqnarray}{r}
\begin{pmatrix} R(\mathbf{i}) \\ G(\mathbf{v}) \end{pmatrix} + \begin{pmatrix} \mathbf{0} & M\tran \\
-M & \mathbf{0}
\end{pmatrix}  \begin{pmatrix} \mathbf{i} \\ \mathbf{v} \end{pmatrix}
 -  \begin{pmatrix} B_R \\ B_G \end{pmatrix}  \mathbf{u} = 0,\label{eq:inclusion1}\\
 \mathbf{y} = (B_R\tran\; B_G\tran)\begin{pmatrix}  \mathbf{i} \\ \mathbf{v} \end{pmatrix}  + D  \mathbf{u},\label{eq:inclusion2}
\end{IEEEeqnarray}
where $M = \tilde{M}\otimes \Id$, $B_R = \tilde{B}_R \otimes \Id$, $B_G = \tilde{B}_G \otimes \Id$ and $D = \tilde{D} \otimes \Id$.
\end{theorem}

It should be pointed out that, while Theorem~\ref{thm:representation} guarantees the existence of a partition of currents and voltages such that the circuit has a monotone+skew form, this partition may not be unique, and for any particular partition, the monotone+skew form of \eqref{eq:inclusion1} and \eqref{eq:inclusion2} is not guaranteed to exist.

Our general approach to solving for the steady state responses will be to first solve Equation~\eqref{eq:inclusion1} for the internal signals $\mathbf{i}$ and $\mathbf{v}$, using a monotone inclusion algorithm, and then to apply the output transformation \eqref{eq:inclusion2} to obtain the response signals in $\mathbf{y}$.

The operators $R$ and $G$ contain impedances and admittances of the circuit elements.  Since we have made no assumption on whether we have an impedance or admittance operator for each element, $R$ and $G$ may contain inverses of operators on $L_{2 ,T}$.  The inverse is taken to be a relational inverse, and always exists, but may not have a domain equal to all of $L_{2, T}$. For example, the admittance of an LTI capacitor is an operator on $L_{2, T}$, but its impedance is a static gain cascaded with an integrator, and only maps signals whose integral on $[0, T]$ is zero into $T$-periodic signals.  The result of this is that, while Equations~\eqref{eq:inclusion1} and~\eqref{eq:inclusion2} can always be written, for a given set of excitation signals, solutions may not exist.  The algorithm we propose in Section~\ref{sec:algorithms} is guaranteed to find a solution, if one exists, but does not guarantee its existence.  It is, however, always well-posed.  Despite the fact that $R$ and $G$ may have restricted domain, they are only ever accessed via their resolvents.  As the elements are assumed to be maximal monotone, it follows from Minty's surjectivity theorem \cite{Minty1961} that their resolvents have domain $L_{2, T}$, and the algorithm steps are always well-defined.

\begin{proof}[Proof of Theorem~\ref{thm:representation}]
    We follow the element extraction approach described in \cite[$\S2$]{Hughes2017c} and illustrated in Figure~\ref{fig:element_extraction}.  We form an $m+n$ circuit $\mathcal{C}_{\text{wires}}$ by replacing each of the elements with an internal port.  Denote the currents at the internal ports by $\{i_k\}_{k = m+1\ldots m+n}$ and the voltages at the internal ports by $\{v_k\}_{k = m+1 \ldots n}$.  The sign of the current at each driving point and internal port is chosen so that the product $i_k(t)v_k(t)$ is the instantaneous power extracted from the port.  It then follows from \cite[Eq. (8)]{Hughes2017c} that there exist $\mathbf{u}, \mathbf{y}, \mathbf{i}, \mathbf{v}$ as specified in the theorem statement and such that, writing $\mathbf{z} = (\mathbf{i}, \mathbf{v})\tran$,
    \begin{IEEEeqnarray}{rCl}
\biggl(\underbrace{\begin{pmatrix} H_{11} & H_{12}\\ H_{21} & H_{22} \end{pmatrix}}_{H}
\otimes \Id\biggr) \begin{pmatrix} \mathbf{u} \\ \mathbf{z} \end{pmatrix}
                       &=& 
                       \begin{pmatrix} \mathbf{y} \\ \tilde{\mathbf{z}} \end{pmatrix},\label{eq:hybrid}
    \end{IEEEeqnarray}
    where $\tilde{\mathbf{z}}$ contains the dual variables to $\mathbf{z}$, that is, if for index $k$ and port $j$, $\mathbf{z}_k$ = $v_j$, then $\tilde{\mathbf{z}}_k = i_j$, and if $\mathbf{z}_k$ = $i_j$, then $\tilde{\mathbf{z}}_k = v_j$.  Since $\mathcal{C}_{\text{wires}}$ is lossless and contains no dynamic components, for any fixed time $t$, $\mathbf{u}\tran(t)\mathbf{y}(t) + \mathbf{z}\tran(t)\tilde{\mathbf{z}}(t) = \sum_{k = 1}^{n + m} i_k(t) v_k(t) = 0$.  It follows that $H$ is skew-symmetric.  Furthermore, because $\mathbf{z}$ is partitioned into $(\mathbf{i}, \mathbf{v})\tran$, $H_{22}$ has the form
    $\left(\begin{smallmatrix} \mathbf{0} & -\tilde{M}\tran \\
    \tilde{M} & \mathbf{0}
    \end{smallmatrix}\right)$,
    where $\tilde{M} \in \R^{q \times p}$ (this follows from the independence of Kirchoff's voltage and current laws, see, for example, \cite[Thm. 6-3]{Seshu1961}).

    The circuit elements give $n$ relations between the entries of $\mathbf{z}$ and the entries of $\tilde{\mathbf{z}}$.  Collecting those that map currents to voltages in a single operator $R: L_{2, T}^p \to  L_{2, T}^p$, and those that map voltages to currents in a single operator $G: L_{2, T}^q \to  L_{2, T}^q$, and substituting $(R(\mathbf{i}), G(\mathbf{v}))\tran$ for $\tilde{\mathbf{z}}$ in Equation~\eqref{eq:hybrid}, gives the equations
    
    \begin{IEEEeqnarray*}{rCl}
\begin{pmatrix} R(\mathbf{i}) \\ G(\mathbf{v}) \end{pmatrix} + \left(\begin{pmatrix} \mathbf{0} & \tilde{M}\tran \\
-\tilde{M} & \mathbf{0}
\end{pmatrix}\otimes \Id\right) \begin{pmatrix} \mathbf{i} \\ \mathbf{v} \end{pmatrix}&&
 \\-  (H_{21} \otimes \Id) \mathbf{u} &=& 0,\\
 \mathbf{y} = (H_{12}\otimes \Id) \begin{pmatrix} \mathbf{i} \\ \mathbf{v} \end{pmatrix} + (H_{11}\otimes \Id) \mathbf{u}.
\label{eq:inclusion}
\end{IEEEeqnarray*}
The proof concludes by noting that skew-symmetry of $H$ implies $H_{12} = -H_{21}\tran$, and writing $\tilde{B} = H_{21}$ and $\tilde{D} = H_{11}$ to obtain Equations~\eqref{eq:inclusion1} and~\eqref{eq:inclusion2}.
\end{proof}

\begin{remark}
Theorem~\ref{thm:representation} shows that the behavior of any monotone $m$-port circuit may be written in the form of \eqref{eq:inclusion1} and \eqref{eq:inclusion2}.  The opposite is also true: the solutions to any set of equations in the form of \eqref{eq:inclusion1} and \eqref{eq:inclusion2} correspond to the behavior of a monotone $m$-port circuit, constructed in the following manner.  Firstly, the matrix $H$ is synthesized using the method of \cite[Ch. 10]{Anderson1973}.  The elements in $R$ and $G$ are then connected across $n$ of the terminal pairs of this realization; the remaining $m$ terminal pairs are the external driving points. 
\end{remark}

\section{MONOTONE INCLUSION ALGORITHMS}\label{sec:algorithms}

Equation~\eqref{eq:inclusion1} expresses the relation between excitations, internal currents and internal voltages as the zero of the sum of two maximal monotone operators, one of which is a structured skew-symmetric matrix.  There are many algorithms for finding the zero of the sum of two maximal monotone operators, including the Douglas-Rachford method and its special case the alternating direction method of multipliers (ADMM), forward-backward splitting (assuming one of the operators is cocoercive), and forward-backward-forward splitting (assuming one of the operators is Lipschitz continuous).  A method that exploits the particular block monotone+skew structure of \eqref{eq:inclusion1} is the Condat--V\~u algorithm~\cite{Condat2013,Vu2013}, which we will use in this paper.  A benefit of this method is that it evaluates $M$ only using forward evaluations, rather than through any sort of inverse.  This gives low-cost iterations that make it suitable for solving large-scale problems.  A second benefit is that the device operators $G$ and $R$ are only evaluated via their resolvents, allowing multi-valued elements such as diodes to be used. This will be further illustrated in Section~\ref{sec:example}.  For our problem \eqref{eq:inclusion1}, the Condat--V\~u iteration \tom{is given by Algorithm~\ref{alg:CV}}.  The algorithm alternately updates the internal currents and internal voltages.
\begin{algorithm}
\caption{The Condat--V\~u algorithm.}
\label{alg:CV}
\begin{algorithmic}
\State \textbf{Given:} zero inclusion of the form \eqref{eq:inclusion1}, step sizes $\sigma, \tau > 0$, initial values $\mathbf{i_0}, \mathbf{v_0}$, tolerance $\epsilon > 0$.\\

\State \textbf{Define:} $\bar{R}(\mathbf{i}) \coloneqq R(\mathbf{i}) - B_R \mathbf{u}$, $\bar{G}(\mathbf{v}) \coloneqq G(\mathbf{v}) - B_G \mathbf{u}$.

\State \textbf{Iterate:}
\begin{align*}
    \mathbf{i}_{k+1}&=(\Id+\tau\bar{R})^{-1}(\mathbf{i}_k - \tau M\tran \mathbf{v}_k)\\
    \mathbf{v}_{k+1}&=(\Id+\sigma\bar{G})^{-1}(\mathbf{v}_k+\sigma M(2\mathbf{i}_{k+1} - \mathbf{i}_k))
\end{align*}
\State \textbf{while} $\|\mathbf{v}_{k+1} - \mathbf{v}_k\| > \epsilon$.
\end{algorithmic}
\end{algorithm}
The fixed points of this iteration are the solutions to \eqref{eq:inclusion1}.  Convergence to a fixed point is guaranteed if $\tau$ and $\sigma$ satisfy $\tau\sigma<\tfrac{1}{\|M\|^2}$, where $\|M\|$ denotes the operator norm of $M$~\cite{Condat2013,Vu2013}.

\section{EXAMPLES}\label{sec:example}

\tom{
In this section, we apply our computational method to two examples: the RLC circuit of Example~\ref{ex:basic}, and a bridge rectifier. Code for both examples can be found at \small{\texttt{https://github.com/ThomasChaffey/circuit-\\analysis-using-monotone-skew-splitting}}.
}

\begin{figure}[ht]
        \centering
        \begin{tikzpicture}
                \begin{axis}
                        [
                        no markers,
                        name = ax1,
                        width=0.45\textwidth,
                        height=0.3\textwidth,
                        ticklabel style={/pgf/number format/fixed},
                        ylabel={V},
                        xlabel={$t$ (s)},
                        cycle list name=colors,
                        grid=both,
                        grid style={line width=.1pt, draw=Gray!20},
                        axis x line=bottom,
                        axis y line=left,
                        scaled x ticks=false,
                        legend pos=south west
                        ]
                        \addplot [dashed] table [x=t, y=vp, col sep = comma, mark =
                                none]{"./basic_example.csv"};
                        \addplot [black] table [x=t, y=vq, col sep = comma, mark =
                                none]{"./basic_example.csv"};
                        \legend{$v_p$, $v_q$}
               \end{axis}
        \end{tikzpicture}
        \caption{\tom{Input voltage $v_p$ and output voltage $v_q$ for the RLC circuit of Figure~\ref{fig:basic_example}, computed using the Condat--V\~u algorithm.  Circuit parameters are $R = 1\, \Omega$, $L = 0.001\,$H and $C = 0.01\,$F, algorithm parameters are $\Delta t = 1\times 10^{-4}$ s (200 samples) and step sizes of $\tau = \sigma = 0.05$.}}%
        \label{fig:basic_example_output}
\end{figure}
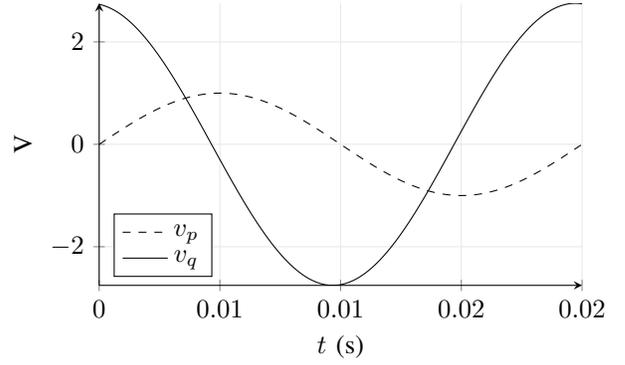

\tom{
\begin{example}\label{ex:basic_continued}
In this example, we revisit the RLC circuit of Example~\ref{ex:basic}, and apply the Condat--V\~u algorithm to compute the voltage response $v_q$ when driving point $p$ is connected to a sinusoidal voltage source, and driving point $q$ is open, so $i_q(t) = 0$ for all $t$. Voltage and current signals are discretized at a regular interval $\Delta t$ to produce discrete signals of length $N$.  The differential operator $\td{}{t}$ is replaced by the periodic backwards finite difference operator, $\nabla(u)(k) \coloneqq (u(k) - u(k-1))/\Delta t$, where $u(-1) \coloneqq u(N)$.  The operators $\hat{C}$ and $\hat{L}$ then have a matrix representation of $CD$ and $LD$ respectively, where $D$ is the $N\times N$ matrix
\begin{IEEEeqnarray*}{rCl}
    \begin{pmatrix}
        1 & 0 & 0 & \ldots & -1 \\
        -1 & 1 & 0 & \ldots & 0 \\
        \vdots & \vdots & \vdots & \ddots & \vdots \\
        0 & 0 & 0 & \ldots & 1
    \end{pmatrix}.
\end{IEEEeqnarray*}
In order to apply the Condat--V\~u algorithm, the resolvents of each element must be calculated.  The resolvents $(\Id + \sigma\hat{C})^{-1}$ and $(\Id + \tau\hat{L})^{-1}$ are given by the matrices $(I + \sigma CD)^{-1}$ and $(I + \tau LD)^{-1}$, which are precomputed and stored in memory.  The resolvent $(\Id + \tau \hat{R})^{-1}$ is equivalent to scalar multiplication by $1/(1 + \tau R)$.  The input offsets $v_p$ and $-i_q$ are incorporated into the resolvents by offsetting their inputs:  $(\Id + \tau\bar{L})^{-1}(i_1) = (\Id + \tau\hat{L})^{-1}(i_1 - \tau v_p)$ and $(\Id + \sigma\bar{R})^{-1}(v_3) = (\Id + \sigma\hat{R})^{-1}(v_3 + \tau i_q)$. To compute $v_q$, given $v_p$ and assuming $i_q(t) = 0$ for all $t$, Algorithm~\ref{alg:CV} is applied to solve Equation~\eqref{eq:inclusion1}, and Equation~\eqref{eq:inclusion2} is then computed.  Example input and output signals are plotted in Figure~\ref{fig:basic_example_output}.
\end{example}
}

\begin{example}\label{ex:rectifier}
    
    This example considers the filtered bridge rectifier shown in Figure~\ref{fig:rectifier}.  We will consider the voltage response $v_q$ when driving point $p$ is connected to a sinusoidal voltage source, and driving point $q$ is connected to a constant current source. We begin by using the element extraction approach of Theorem~\ref{thm:representation} to express the circuit in monotone+skew form.

\begin{figure}
\centering
    \includegraphics{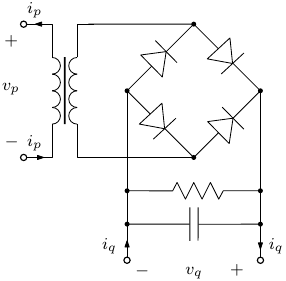}
    \caption{A bridge rectifier with filtering capacitor (capacitance $C$) and nominal load resistor (resistance $R$).  The transformer is a $1:24$ step-down transformer.}
    \label{fig:rectifier}
\end{figure}

\begin{figure}
\centering
    \includegraphics{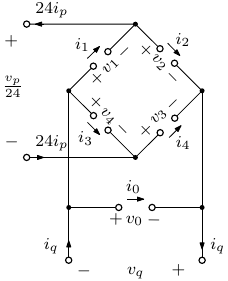}
    \caption{The bridge rectifier of Figure~\ref{fig:rectifier} with elements replaced by internal ports.  Polarities at the ports are chosen so that the product $i_k v_k$ represents power extracted from the circuit.}
    \label{fig:rectifier_wires}
\end{figure}

Figure~\ref{fig:rectifier_wires} shows the circuit of Figure~\ref{fig:rectifier} with the elements replaced by internal ports.  The resistor and capacitor are combined into a single subcircuit, and replaced by port~0.  The transformer results in extra factors of $24$ at driving point $p$.  Following the procedure of Section~\ref{sec:representation}, we select an excitation and response variable at each driving point, and derive a hybrid representation for the circuit of wires.  $v_p$ is selected as the excitation at driving point~$p$ and $v_q$ is selected as the response at driving point~$q$.  The other variables may be partitioned arbitrarily, provided a hybrid representation exists for the partition.  For example, we can choose the currents at ports~0, 1 and~2, and the voltages at ports~3 and~4, to be excitation variables.  This gives the hybrid representation

\begin{IEEEeqnarray}{rCl}
\begin{pmatrix} i_p \\ v_q \\ v_0 \\ v_1 \\ v_2 \\ i_3 \\ i_4 \end{pmatrix} &{=}&
\left(\begin{array}{c c | c c c c c} 0 & 0 & 0 & \frac{1}{24} & \frac{-1}{24} & 0 & 0 \\
                0 & 0 & 0 & 0 & 0 & -1 & -1 \\\hline
                0 & 0 & 0 & 0 & 0 & -1 & -1 \\
                \frac{-1}{24} & 0 & 0 & 0 & 0 & 0 & 1 \\
                \frac{1}{24} & 0 & 0 & 0 & 0 & 1 & 0 \\
                0 & 1 & 1 & 0 & -1 & 0 & 0 \\
                0 & 1 & 1 & -1 & 0 & 0 & 0 
\end{array}\right)
\begin{pmatrix} v_p \\ i_q \\ i_0 \\ i_1 \\ i_2 \\ v_3 \\ v_4 \end{pmatrix}
\label{eq:wires}
\end{IEEEeqnarray}

where the time dependence of each signal has been dropped.  The devices are described by
\begin{IEEEeqnarray}{rCl}
\begin{pmatrix} v_0 \\ v_1 \\ v_2 \\ i_3 \\ i_4 \end{pmatrix} &=& \begin{pmatrix} R(\mathbf{i}) \\ G(\mathbf{v})\end{pmatrix} \coloneqq
\begin{pmatrix} R_{RC} (i_0) \\ R_D(i_1) \\ R_D(i_2) \\ G_D(v_3) \\ G_D(v_4) \end{pmatrix},
\label{eq:devices}
\end{IEEEeqnarray}
where $\mathbf{i} \coloneqq (i_0, i_1, i_2)\tran$ and $\mathbf{v} \coloneqq (v_3, v_4)\tran$.
$R_{RC}$ is the convolution operator of the RC filter in impedance form, described by the operator $R(RC\td{}{t} + 1)^{-1}$, where $\td{}{t}$ represents the differential operator on $L_{2, T}$. $R_D$ is the impedance of an ideal diode:
\begin{IEEEeqnarray*}{rCl}
    (R_D i)(t) &=& \begin{cases} \{0\} & i(t) > 0\\
                        (-\infty, 0] & i(t) = 0.
                           \end{cases}
\end{IEEEeqnarray*}
$G_D$ is the admittance of an ideal diode, given by the relational inverse of $R_D$:
\begin{IEEEeqnarray*}{rCl}
    (G_D v)(t) &=& \begin{cases} \{0\} & v(t) < 0\\
                           [0, \infty) & v(t) = 0.
                           \end{cases}
\end{IEEEeqnarray*}
Substituting Equation~\eqref{eq:devices} in Equation~\eqref{eq:wires} gives
\begin{IEEEeqnarray}{rCl}
\begin{pmatrix} R(\mathbf{i}) \\ G(\mathbf{v}) \end{pmatrix} + \begin{pmatrix} \mathbf{0}_{3\times 3} & M\tran \\
-M & \mathbf{0}_{2\times 2}
\end{pmatrix} \begin{pmatrix} \mathbf{i} \\ \mathbf{v} \end{pmatrix}
 -  \begin{pmatrix} B_R \\ B_G \end{pmatrix} \mathbf{u} = 0,
\label{eq:example_inclusion}
\end{IEEEeqnarray}
where
\begin{IEEEeqnarray*}{l}
\mathbf{u} \coloneqq \begin{pmatrix} v_p \\ i_q \end{pmatrix}, \qquad
M \coloneqq \begin{pmatrix} 
                 1 & 0 & -1  \\
                 1 & -1 & 0  
\end{pmatrix}\otimes \Id,\\
B_R \coloneqq \begin{pmatrix} 0 & 0 \\ \frac{-1}{24} & 0 \\ \frac{1}{24} & 0 \end{pmatrix}\otimes \Id, \qquad B_G \coloneqq \begin{pmatrix} 0 & 1 \\ 0 & 1 \end{pmatrix}\otimes \Id.
\end{IEEEeqnarray*}
$v_q$ and $i_p$ are then obtained by the output equations $v_q = -v_0 = -v_3 - v_4$ and $i_p = i_1 - i_2$.

Similarly to Example~\ref{ex:basic_continued}, the voltage and current signals are discretized at a regular interval $\Delta t$ to produce discrete signals of length $N$, and the differential operator $\td{}{t}$ is replaced by the periodic backwards finite difference operator.
In order to apply the Condat--V\~u algorithm, the resolvents $(\Id + \tau \bar{R} )^{-1}$ and $(\Id + \sigma \bar{G} )^{-1}$ must be calculated.  It follows from \cite[Prop. 23.16]{Bauschke2011} that the resolvents can be calculated elementwise.  Specifically, we have:
\begin{IEEEeqnarray*}{rCl}
\begin{pmatrix} (\Id + \tau R )^{-1} \\ (\Id + \sigma G )^{-1} \end{pmatrix}  =
\begin{pmatrix} (\Id + \tau R_{RC})^{-1} \\ (\Id + \tau R_D)^{-1} \\ (\Id + \tau R_D)^{-1} \\ (\Id + \sigma G_D)^{-1} \\ (\Id + \sigma G_D)^{-1}\end{pmatrix}.
\end{IEEEeqnarray*}
The offsets are dealt with by offsetting the inputs to the resolvent operators: $(\Id + \tau \bar{R} )^{-1}(\mathbf{i}) =(\Id + \tau R )^{-1}(\mathbf{i} + B_R \mathbf{u})$, and similarly for $\bar{G}$.

The resolvent of the $RC$ circuit impedance has an explicit form as multiplication by a matrix:
$(\Id + \tau R_{RC})^{-1}(i) = (I + \tau R(RC\nabla + 1)^{-1})^{-1} i$.  The resolvent of the impedance of a diode is a rectified linear unit (ReLU): 
\begin{IEEEeqnarray*}{rCl}
(\Id + \tau R_D)^{-1}(i)(t) = \operatorname{ReLU}(i(t)) \coloneqq \begin{cases} i(t) & i(t) > 0\\
                                            0 & \text{otherwise}.
                                            \end{cases}
\end{IEEEeqnarray*}
This can be most easily seen graphically by plotting the current-voltage characteristic (Figure~\ref{fig:resolvent_diode}).  Similarly, the resolvent of the admittance of a diode is given by $(\Id + \sigma G_D)^{-1}(v)(t) = -\operatorname{ReLU}(-v(t))$.

\begin{figure}[ht]
    \centering
    \includegraphics{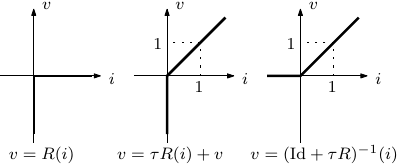}
    \caption{The $\tau$-resolvent of the impedance of a diode is a rectified linear unit (ReLU).}
    \label{fig:resolvent_diode}
\end{figure}

Applying the Condat--V\~u algorithm, with a resistance of $R = 1\times 10^3\, \Omega$, a capacitance of $C = 10\,\mu$F, a 50 Hz sinuosidal input voltage $v_p(t) = 240.0 \sin(100\pi)$ V, a constant input current $i_q(t) = -5$ mA, a discretization interval of $\Delta t = 1\times 10^{-4}$ s (200 samples) and step sizes of $\tau = \sigma = 0.005$ gives the output shown in Figures~\ref{fig:result1} and~\ref{fig:result2}.  
\end{example}

\begin{figure}[ht]
        \centering
        \begin{tikzpicture}
                \begin{axis}
                        [
                        no markers,
                        name = ax1,
                        width=0.45\textwidth,
                        height=0.3\textwidth,
                        ticklabel style={/pgf/number format/fixed},
                        ylabel={V},
                        xlabel={$t$ (s)},
                        cycle list name=colors,
                        grid=both,
                        grid style={line width=.1pt, draw=Gray!20},
                        axis x line=bottom,
                        axis y line=left,
                        scaled x ticks=false,
                        legend pos=south west
                        ]
                        \addplot [dashed] table [x=t, y=vp, col sep = comma, mark =
                                none]{"./bridge_rectifier.csv"};
                        \addplot [black] table [x=t, y=vq, col sep = comma, mark =
                                none]{"./bridge_rectifier.csv"};
                        \legend{$v_p/24$, $v_q$}
               \end{axis}
        \end{tikzpicture}
        \caption{Input voltage $v_p$ and output voltage $v_q$ for the filtered bridge rectifier of Figure~\ref{fig:rectifier}, computed using the Condat--V\~u algorithm.  Circuit parameters are $R = 1\times 10^3\, \Omega$ and $C = 10\,\mu$F, algorithm parameters are $\Delta t = 1\times 10^{-4}$ s (200 samples) and step sizes of $\tau = \sigma = 0.005$.}%
        \label{fig:result1}
\end{figure}
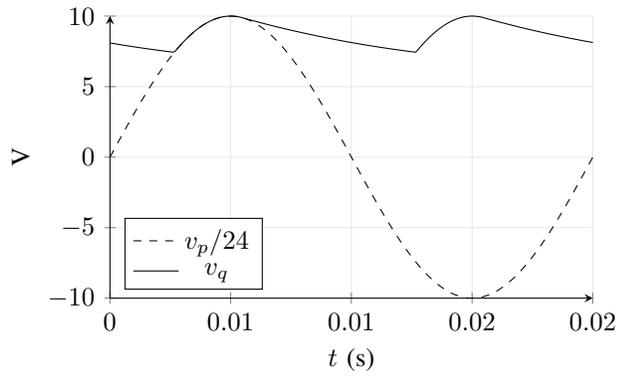

\begin{figure}[ht]
        \centering
        \begin{tikzpicture}
                \begin{axis}
                        [
                        no markers,
                        name = ax1,
                        width=0.45\textwidth,
                        height=0.3\textwidth,
                        ticklabel style={/pgf/number format/fixed},
                        ylabel={A},
                        xlabel={$t$ (s)},
                        cycle list name=colors,
                        grid=both,
                        grid style={line width=.1pt, draw=Gray!20},
                        axis x line=bottom,
                        axis y line=left,
                        scaled x ticks=false,
                        legend pos=south east
                        ]
                        \addplot [dashed] table [x=t, y=iq, col sep = comma, mark =
                                none]{"./bridge_rectifier.csv"};
                        \addplot [black] table [x=t, y=ip, col sep = comma, mark =
                                none]{"./bridge_rectifier.csv"};
                        \legend{$i_q$, $24i_p$}
               \end{axis}
        \end{tikzpicture}
        \caption{Input current $i_q$ and output current $i_p$ for the filtered bridge rectifier of Figure~\ref{fig:rectifier}.  Algorithm and circuit parameters as for Figure~\ref{fig:result1}.}%
        \label{fig:result2}
\end{figure}
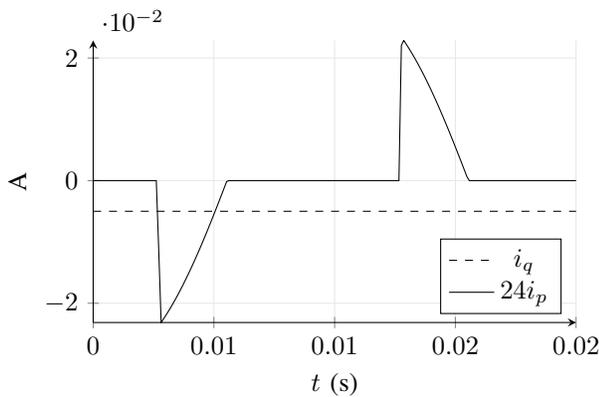

\section{CONCLUSIONS}
 
 We have demonstrated that the behavior of a nonlinear $m$-port circuit composed of monotone elements may be expressed as the zero of the sum of a structured skew linear operator and a monotone operator.  This is precisely the form of zero inclusion which may be solved using the Condat--V\~u algorithm.  We have used this correspondence to develop a method for solving the periodic forced behavior of a monotone $m$-port circuit, extending the method of \cite{Chaffey2021d} to arbitrary numbers of ports and arbitrary interconnection structures.
 
 The method of \cite{Chaffey2021d} has been adapted to solve for unforced periodic solutions in one-port circuits composed of the \emph{difference} of monotone elements by Das, Chaffey and Sepulchre \cite{Das2022}.  Such circuits include the van der Pol and FitzHugh--Nagumo oscillators.  An interesting avenue for future work is to extend the the skew+monotone method we present here in a similar manner, to differences of monotone $m$-port circuits.  A second avenue for future research is the extension to circuits containing active components, such as ideal operational amplifiers.

 \tom{A final question for future is whether the monotone+skew form of a nonlinear circuit allows for algorithmic solutions to problems other than operator inversion or simulation. For example, are there optimal control problems that can be solved using a splitting method similar to that used in this paper?}


\bibliographystyle{ieeetr}
\bibliography{skew-monotone-bibtex.bib}

\begin{thebibliography}{10}

\bibitem{Duffin1946}
R.~J. Duffin, ``Nonlinear networks. {{I}},'' {\em Bulletin of the American
  Mathematical Society}, vol.~52, no.~10, pp.~833--839, 1946.

\bibitem{Minty1960}
G.~J. Minty, ``Monotone networks,'' {\em Proceedings of the Royal Society of
  London. Series A. Mathematical and Physical Sciences}, vol.~257, no.~1289,
  pp.~194--212, 1960.

\bibitem{Minty1961}
G.~J. Minty, ``Solving steady-state nonlinear networks of ``monotone''
  elements,'' {\em IRE Transactions on Circuit Theory}, vol.~8, no.~2,
  pp.~99--104, 1961.

\bibitem{Rockafellar1967}
R.~T. Rockafellar, ``Convex programming and systems of elementary monotonic
  relations,'' {\em Journal of Mathematical Analysis and Applications},
  vol.~19, no.~3, pp.~543--564, 1967.

\bibitem{Rockafellar1976}
R.~T. Rockafellar, ``Monotone operators and the proximal point algorithm,''
  {\em SIAM Journal on Control and Optimization}, vol.~14, no.~5, pp.~877--898,
  1976.

\bibitem{Rockafellar1997}
R.~T. Rockafellar and R.~J.-B. Wets, {\em Variational {{Analysis}}}.
\newblock Springer, 1997.

\bibitem{Ryu2016}
E.~K. Ryu and S.~Boyd, ``A primer on monotone operator methods,'' {\em Appl.
  Comput. Math.}, vol.~15, no.~1, pp.~3--43, 2016.

\bibitem{Ryu2022}
E.~Ryu and W.~Yin, {\em Large-Scale Convex Optimization: {{Algorithms}} \&
  Analyses via Monotone Operators}.
\newblock {Cambridge University Press}.

\bibitem{Parikh2013}
N.~Parikh and S.~Boyd, ``Proximal algorithms,'' {\em Foundations and Trends in
  Optimization}, vol.~1, no.~3, pp.~123--231, 2013.

\bibitem{Bertsekas2011}
D.~P. Bertsekas, ``Incremental proximal methods for large scale convex
  optimization,'' {\em Mathematical Programming}, vol.~129, no.~2,
  pp.~163--195, 2011.

\bibitem{Combettes2011}
P.~L. Combettes and J.-C. Pesquet, ``Proximal splitting methods in signal
  processing,'' in {\em Fixed-Point Algorithms for Inverse Problems in Science
  and Engineering}, Springer {{Optimization}} and {{Its Applications}},
  pp.~185--212, {New York}: {Springer}, 2011.

\bibitem{Brogliato2020}
B.~Brogliato and A.~Tanwani, ``Dynamical systems coupled with monotone
  set-valued operators: Formalisms, applications, well-posedness, and
  stability,'' {\em SIAM Review}, vol.~62, no.~1, pp.~3--129, 2020.

\bibitem{Camlibel2016}
M.~K. \c{C}amlibel and J.~M. Schumacher, ``Linear passive systems and maximal
  monotone mappings,'' {\em Mathematical Programming}, vol.~157, no.~2,
  pp.~397--420, 2016.

\bibitem{Brogliato2004}
B.~Brogliato, ``Absolute stability and the {{Lagrange\textendash Dirichlet}}
  theorem with monotone multivalued mappings,'' {\em Systems \& Control
  Letters}, vol.~51, no.~5, pp.~343--353, 2004.

\bibitem{Brezis1973}
H.~Brézis, {\em Opérateurs maximaux monotones et semi-groupes de contractions
  dans les espaces de {{Hilbert}}}.
\newblock No.~5 in North-Holland Mathematics Studies, Amsterdam: North-Holland
  Pub. Co, 1973.

\bibitem{Burger2014}
M.~Bürger, D.~Zelazo, and F.~Allgöwer, ``Duality and network theory in
  passivity-based cooperative control,'' {\em Automatica}, vol.~50, no.~8,
  pp.~2051--2061, 2014.

\bibitem{Chaffey2021d}
T.~Chaffey and R.~Sepulchre, ``Monotone one-port circuits,'' {\em IEEE
  Transactions on Automatic Control}, 2023 (early access).

\bibitem{Krylov1947}
N.~S. Krylov and N.~N. Bogoliubov, {\em Introduction to Non-Linear Mechanics}.
\newblock {Princeton}: {Princeton University Press}, 1947.

\bibitem{Aprille1972}
T.~Aprille and T.~Trick, ``Steady-state analysis of nonlinear circuits with
  periodic inputs,'' {\em Proceedings of the IEEE}, vol.~60, no.~1,
  pp.~108--114, 1972.

\bibitem{Cellier2006}
F.~E. Cellier and E.~Kofman, {\em Continuous System Simulation}.
\newblock {New York}: {Springer}, 2006.

\bibitem{Acary2008}
V.~Acary and B.~Brogliato, {\em Numerical Methods for Nonsmooth Dynamical
  Systems: Applications in Mechanics and Electronics}, vol.~35 of {\em Lecture
  {{Notes}} in {{Applied}} and {{Computational Mechanics}}}.
\newblock Springer, 2008.

\bibitem{Heemels2017}
W.~Heemels, V.~Sessa, F.~Vasca, and M.~\c{C}amlibel, ``Computation of periodic
  solutions in maximal monotone dynamical systems with guaranteed
  consistency,'' {\em Nonlinear Analysis: Hybrid Systems}, vol.~24,
  pp.~100--114, 2017.

\bibitem{Iannelli2011}
L.~Iannelli, F.~Vasca, and G.~Angelone, ``Computation of steady-state
  oscillations in power converters through complementarity,'' {\em IEEE
  Transactions on Circuits and Systems I: Regular Papers}, vol.~58, no.~6,
  pp.~1421--1432, 2011.

\bibitem{Meingast2014}
M.~B. Meingast, M.~Legrand, and C.~Pierre, ``A linear complementarity problem
  formulation for periodic solutions to unilateral contact problems,'' {\em
  International Journal of Non-Linear Mechanics}, vol.~66, pp.~18--27, 2014.

\bibitem{vanderSchaft2014}
A.~{van der Schaft} and D.~Jeltsema, {\em Port-Hamiltonian Systems Theory: An
  Introductory Overview}.
\newblock Foundations and Trends in Systems and Control, Now Foundations and
  Trends, 2014.

\bibitem{Briceno-Arias2010}
L.~{Brice\~{n}o-Arias} and P.~L. Combettes, ``A monotone+skew splitting model
  for composite monotone inclusions in duality,'' 2010.

\bibitem{Condat2013}
L.~Condat, ``A primal\textendash dual splitting method for convex optimization
  involving lipschitzian, proximable and linear composite terms,'' {\em Journal
  of Optimization Theory and Applications}, vol.~158, no.~2, pp.~460--479,
  2013.

\bibitem{Vu2013}
B.~C. V{\~u}, ``A splitting algorithm for dual monotone inclusions involving
  cocoercive operators,'' {\em Advances in Computational Mathematics}, vol.~38,
  no.~3, pp.~667--681, 2013.

\bibitem{Chambolle2011}
A.~Chambolle and T.~Pock, ``A first-order primal-dual algorithm for convex
  problems with applications to imaging,'' {\em Journal of Mathematical Imaging
  and Vision}, vol.~40, no.~1, pp.~120--145, 2011.

\bibitem{Desoer1975}
C.~A. Desoer and M.~Vidyasagar, {\em Feedback Systems: Input\textendash Output
  Properties}.
\newblock {Elsevier}, 1975.

\bibitem{Willems2007}
J.~C. Willems, ``The behavioral approach to open and interconnected systems,''
  {\em IEEE Control Systems}, vol.~27, no.~6, pp.~46--99, 2007.

\bibitem{Desoer1969}
C.~A. Desoer and E.~S. Kuh, {\em Basic Circuit Theory}.
\newblock {McGraw-Hill}.

\bibitem{Hughes2017c}
T.~H. Hughes, ``Passivity and electric circuits: A behavioral approach,'' {\em
  IFAC-PapersOnLine}, vol.~50, no.~1, pp.~15500--15505, 2017.

\bibitem{Anderson1966}
B.~Anderson and R.~Newcomb, ``Cascade connection for time-invariant $n$-port
  networks,'' {\em Proceedings of the Institution of Electrical Engineers},
  vol.~113, no.~6, p.~970, 1966.

\bibitem{Seshu1961}
S.~Seshu and M.~B. Reed, {\em Linear Graphs and Electrical Networks}.
\newblock {Addison-Wesley}, 1961.

\bibitem{Anderson1973}
B.~D.~O. Anderson and S.~Vongpanitlerd, {\em Network Analysis and Synthesis: A
  Modern Systems Theory Approach}.
\newblock Prentice-Hall Networks Series, {Englewood Cliffs, N.J.}:
  {Prentice-Hall}, 1973.

\bibitem{Bauschke2011}
H.~H. Bauschke and P.~L. Combettes, {\em Convex Analysis and Monotone Operator
  Theory in Hilbert Spaces}.
\newblock {{CMS Books}} in {{Mathematics}}, {New York, NY}: {Springer New
  York}, 2011.

\bibitem{Das2022}
A.~Das, T.~Chaffey, and R.~Sepulchre, ``Oscillations in mixed-feedback
  systems,'' {\em Systems \& Control Letters}, vol.~166, 2022.

\end{thebibliography}
\end{document}